\documentclass[12pt]{amsart} 

\usepackage{amsmath,amssymb,amscd,amsthm,fullpage, graphicx}
\usepackage{latexsym}
\usepackage{colortbl}

\newtheorem{theorem}{Theorem}

\newtheorem{definition}{Definition}
\newtheorem{remark}{Remark}
\newtheorem{corollary}{Corollary}

\newtheorem{lemma}{Lemma}

\newcommand{\beginproof}{\noindent{\bf Proof: }}

\newcommand{\vv}{{\boldsymbol v}}
\newcommand{\vz}{{\boldsymbol z}}
\newcommand{\vx}{{\boldsymbol x}}
\newcommand{\ve}{{\boldsymbol e}}
\newcommand{\vu}{{\boldsymbol u}}

\newcommand{\vxi}{{\boldsymbol \xi}}
\newcommand{\vy}{{\boldsymbol y}}
\newcommand{\vn}{{\boldsymbol n}}
\newcommand{\vgamma}{{\boldsymbol \gamma}}

\newcommand{\SP}{{\mathbb S}}

\newcommand{\LLL}{{\mathcal L}}

\newcommand{\card}{\mbox{card}}
\newcommand{\spn}{\text{span}}
\def\RR{{\mathbb R}}
\def\N{{\mathbb N}}

\def\ZZ{{\mathbb Z}}

\def\G{{\rm G}}

\def\vol{\mbox{\rm vol}}

\def\spn{\hbox{\rm span}\,}

\def\conv{\hbox{\rm conv}}

\def\endproof{\begin{flushright}
$ \Box $ \\
\end{flushright}}

\begin{document}
\hfill\today
\bigskip
\author{Matthew Alexander}\thanks{ The first author is supported in part by the U.S. National Science
Foundation Grant DMS-1101636.} 
\author{Martin Henk}  
\author{Artem Zvavitch}\thanks{The third author is supported in part by the U.S. National Science Foundation Grant DMS-1101636 and by the  Simons Foundation.}

\address{Department of Mathematics, Kent State University,
Kent, OH 44242, USA} \email{malexan5@kent.edu}

\address{Technische Universit\"at Berlin, Institut f\"ur Mathematik, Strasse des 17.
Juni 136, D-10623 Berlin, Germany} \email{henk@math.tu-berlin.de}

\address{Department of Mathematics, Kent State University,
Kent, OH 44242, USA} \email{zvavitch@kent.edu}

\title{A discrete  version of Koldobsky's slicing inequality}

\begin{abstract}
Let $\# K$ be a number of integer lattice points  contained in a set $K$. In this paper we prove that for each $d\in \N$ there exists a constant $C(d)$ depending on  $d$ only, such that for any origin-symmetric convex body $K \subset \RR^d$ containing $d$ linearly independent lattice points  
$$
 \# K \leq C(d)\max(\# (K\cap H))\, \vol_d(K)^{\frac{d-m}{d}},
$$
where the maximum is taken over all $m$-dimensional  subspaces of $\RR^d$.
We also prove that  $C(d)$ can be chosen asymptotically of order $O(1)^{d}d^{d-m}$. In
addition, we show that if  $K$ is an unconditional convex body then $C(d)$ can be chosen asymptotically of order $O(d)^{d-m}$. 

\end{abstract}

\maketitle

\section{Introduction}

As usual, we will say that $K\subset \RR^d$ is a convex body if $K$ is a convex, compact subset of $\RR^d$ equal to the closure of its interior.  We say that $K$ is origin-symmetric if $K=-K$, where $\lambda K=\{\lambda \vx: \vx \in K\},$ for $\lambda \in \RR$. For a set $K$
we denote by dim$(K)$ its dimension, that is, the dimension of the affine hull of $K$.  We define  $K+L=\{\vx+\vy: \vx\in K, \vy \in L\}$  to be the Minkowski sum of  $K, L \subset \RR^d$.  We will also denote by $\vol_d$ the $d$-dimensional Hausdorff measure, and if the body $K$ is $d$-dimensional we will call $\vol_d(K)$ the  volume of $K$. Finally, let us denote by $\vxi^\perp$ a hyperplane perpendicular to a unit vector $\vxi$, i.e.
$$
\vxi^\perp=\{\vx\in \RR^d:  \vx\cdot \vxi =0\}.
$$
 We refer to \cite{Ga, K3, BGVV, RZ, Sch} for general definitions and properties of convex bodies.

The slicing problem  of Bourgain \cite{Bo1, Bo2} is, undoubtedly, one of the major open problems in convex geometry asking if a convex, origin-symmetric body of volume one must have a large (in volume) hyperplane section. More precisely, it asks whether there exists an absolute constant $\LLL_1$ so that for any origin-symmetric convex body $K$ in $\RR^d$
\begin{equation}\label{eq:kold}
\vol_d(K)^{\frac {d-1}{d}} \le \LLL_1 \max_{\vxi \in \SP^{d-1}} \vol_{d-1}(K\cap \vxi^\bot).
\end{equation}
 The  problem is still open, with the best-to-date estimate of $\LLL_1\le O(d^{1/4})$ established
by Klartag \cite{Kl}, who improved the previous estimate of Bourgain \cite{Bo2}, we refer to  \cite{MP} and  \cite{BGVV} for detailed information and history of the problem.  Recently, Alexander Koldobsky proposed an interesting generalization of the  slicing problem \cite{K1,K2,K4,K5,K6}:
Does there exists an absolute constant $\LLL_2$ so that for every even measure $\mu$
on $\RR^d$, with a positive density, and for every origin-symmetric convex body $K$ in $\RR^d$ such that
\begin{equation}
\label{eq:cont}\mu(K)\le \LLL_2\max_{\vxi\in \SP^{d-1}} \mu(K\cap \vxi^\bot) \vol_d(K)^{\frac 1d}?
\end{equation}
Koldobsky was able to solve the above question for a number of special cases of the body $K$ and provide  a general estimate of $O(\sqrt{d})$. The most amazing fact here is that the constant $\LLL_2$ in (\ref{eq:cont}) can be chosen independent of the measure $\mu$ under the assumption that $\mu$ has even positive density. In addition,  Koldobsky and the second named author were able to prove in \cite{KoZ} that $\LLL_2$ is of order $O\left(d^{1/4}\right)$  if one assumes that the measure $\mu$ is $s$-concave.  We note that the assumption of positive density is essential for the above results and  (\ref{eq:cont}) is simply not true if this condition is dropped. Indeed, to create a counterexample consider an even measure $\mu$ on $\RR^2$ uniformly distributed over $2N$ points on the unit circle, then the constant $\LLL_2$ in (\ref{eq:cont}) will depend on $N$.

During the 2013 AIM workshop  on ``Sections of convex bodies'' Koldobsky asked if  it is possible to provide a discrete analog of inequality (\ref{eq:cont}): Let $\ZZ^d$ be the standard integer lattice in $\RR^d$, $K$ be a convex, origin-symmetric body,  define $\#K=\card(K\cap\ZZ^d)$, the number of points of $\ZZ^d$ in $K$.
\vskip  1em
\noindent{\bf Question:} {\it Does there exist a  constant $\LLL_3$ such that}
$$
\#K\leq \LLL_3 \max_{\vxi\in \SP^{d-1}} \left( \#(K\cap \vxi^\perp)\right)  \vol_d(K)^{\frac{1}{d}},
$$
{\it for all convex origin-symmetric bodies  $K\subset \RR^d$
  containing  $d$ linearly independent lattice points?}
\vskip  1em
We note here that we require that $K$ contains $d$ linearly
independent lattice points, i.e., $\dim(K\cap\ZZ^d)=d$, in order to eliminate the
degenerate case of a body (for example, take a box $[-1/n,
1/n]^{d-1}\times [-20, 20]$) whose maximal section contains all
lattice points in the body, but whose volume may be taken to 0 by eliminating a dimension.

Koldobsky's question is yet another example of  an attempt to translate questions and facts from classical Convexity to more general settings including Discrete Geometry.  The properties of sections of convex bodies with respect to the integer lattice were extensively studied in Discrete Tomography   \cite{GGroZ, GGr1, GGr2, GGro},  where many interesting new properties were proved and a series of exciting open questions were proposed. It is interesting to note that after translation many questions become quite non-trivial and counterintuitive, and the answer may be quite different from the continuous case. In addition, finding the relation between the geometry of a convex set and the number of integer points contained in the set is always a non-trivial task. One can see this, for example, from the history of Khinchin's flatness theorem \cite{Ba1, Ba2, BLPS, KL}.

The main goal of this paper is to study Koldobsky's question.
In Section 2  we will show the solution for the $2$-dimensional
case. The solution is based on the classical Minkowski's First and
Pick's theorems from the Geometry of Numbers and gives a general idea
of the approach to be used in Sections 3 and 4. In Section 3  we apply
a discrete  version of the theorem of F. John  due to T. Tao and V. Vu
\cite{TV1} to give a partial answer to Koldobsky's question and show
that the constant $\LLL_3$ can be chosen independent of the body $K$ and
as small as  $O(d)^{7d/2}$.  We start Section 4  with a case of
unconditional bodies and present a simple proof that in this case
$\LLL_3$ can be chosen of order  $O(d)$ which is  best possible. After, we prove the
discrete analog of Brunn's theorem and use it to show that the constant
$\LLL_3$, for the general case, can be chosen as small as  $O(1)^d$.  
In fact, we prove the  slightly more general result that 
$$
 \# K \leq O(1)^d d^{d-m} \max\left(\# (K\cap H)\right)\,\, \vol_d(K)^{\frac{d-m}{d}},
$$
where the maximum is taken over all $m$-dimensional linear  subspaces $H \subset \RR^d$.   
Finally,  we also provide a short  observation that $\LLL_1 \le \LLL_3$. 

\smallskip 
\noindent {\bf Acknowledgment}: We are indebted to Alexander Koldobsky and  Fedor Nazarov  for
valuable discussions.

\section{Solution in $\ZZ^2$}

Let us start with recalling two classical statements in the Geometry
of Numbers 
(see \cite{TV2}, Theorem 3.28 pg 134 and \cite{BR}, Theorem 2.8 pg 38):

\begin{theorem}\label{th:M1} (Minkowski's First Theorem)
Let  $K \subset \RR^d$ be an origin-symmetric convex body such that $\vol_d(K)\ge 2^d$ then $K$ contains at least one non-zero element of $\ZZ^d$.
\end{theorem}

\begin{theorem} (Pick's Theorem)
Let $P$ be an integral $2$-dimensional convex polygon, then $A=I+\frac{1}{2}B-1$ where $A=\vol_2(P)$ is the area of the polygon, $I$ is the number of lattice points in the interior of $P$, and $B$ is the number of lattice points on the boundary.
\end{theorem}
\noindent  Here a polygon is called integral if it can be described as  the
convex hull of lattice points. 

Now we will use the above theorems to show that the constant $\LLL_3$ in Koldobsky's question can be chosen independently of a convex, origin-symmetric body in $\RR^2$.

\begin{theorem} Let $K$ be a convex origin-symmetric body in $\RR^2$,
  $\dim(K\cap\ZZ^2)=2$, then
$$
\#K\leq 4 \max_{\vxi\in \SP^{1}}  \#(K\cap \vxi^\perp)\,\vol_2(K)^{\frac{1}{2}}.
$$
\end{theorem}
\beginproof Let $s=\sqrt{\vol_2(K)/4}$, then by Minkowski's theorem, since $\vol_2(\frac{1}{s}K)=4$, there exists a non-zero vector $\vu\in \ZZ^2 \cap \frac{1}{s} K$.
Then $s\vu \in K$ and
$$\# \left(L_\vu \cap K\right) \geq 2 \lfloor s \rfloor + 1 
,$$ where $\lfloor s \rfloor$ is the integer part of $s$, 
and $L_\vu$ is the line containing $\vu$ and the origin.
Next, consider $P=\conv(K\cap \ZZ^2)$, i.e., the convex hull of the
integral points inside $K$.   $P$ is an integral $2$-dimensional  convex polytope, and so by Pick's theorem we get that
$$
\vol_2(P)=I+\frac{1}{2}B-1\geq \frac{I+B}{2}-\frac{1}{2},
$$
using that $I \ge 1$.  Thus
$$
\# P=I+B   \leq 2\vol_2(P)+1\leq \frac{5}{2}\vol_2(P),$$
since the minmal volume of an origin-symmetric
integral convex polygon is at least 2.   
We now have that 
\begin{equation*}
\begin{split} 
\# K & = \# P \leq \frac{5}{2}\vol_2(P) \leq \frac{5}{2}\vol_2(K)\\ 
&\leq \frac{5}{2} \,(2\,s)\, \vol_2(K)^\frac{1}{2} < 4\,( 2 \lfloor s \rfloor + 1)\,
\vol(K)_2^\frac{1}{2}\\ 
&\leq  4 \max_{\vxi\in \SP^{1}}  \#(K\cap \vxi^\perp)\,\vol_2(K)^{\frac{1}{2}}.
\end{split}
\end{equation*}
\endproof
	
\section{Approach via Discrete F. John Theorem}

It is a standard technique to get a first estimate in slicing inequalities, i.e. $\LLL_1 \le O(\sqrt{d})$,  by using the classical F. John theorem, \cite{J}, \cite{MS}, or \cite{BGVV},  which claims that for every convex origin-symmetric body $K\subset \RR^d$ there exists an Ellipsoid $E$ such that
$ E\subset K \subset \sqrt{d} E$. 
In this section we will use a recent discrete version of F. John's theorem, proved by T. Tao and V. Vu (see \cite{TV1, TV2}) to prove that the constant $\LLL_3$ in Koldobsky's question  can be chosen independent  of the origin-symmetric convex body $K\subset \RR^d$. We first recall the definition of a generalized arithmetic progression (see \cite{TV1, TV2} for more details):

\begin{definition}
Let $G$ be an additive group, $N=(N_1, \ldots, N_d)$ an $d$-tuple of non-negative integers and $\vv=(\vv_1, \ldots, \vv_d)\in G^d$.
Then a generalized symmetric arithmetic progression ${\bf P}$ is a triplet $(N, \vv, d)$.  In addition, define
$$\text{\rm Image}({\bf P})=  [-N, N] \cdot \vv = \left\{n_1 \vv_1 + \ldots + n_d \vv_d : n_j \in [-N_j, N_j] \cap \ZZ  \text{ for all } 1\leq j \leq d \right\}.$$
The progression is called proper if the map $\vn \mapsto \vn\cdot \vv$
is injective, $\vv=(\vv_1, \ldots, \vv_d)$ is called  its basis
vectors, and $d$ its rank.
\label{def:progression}
\end{definition}

Below is a version for $\ZZ^d$ of the Discrete John theorem from \cite{TV1} (Theorem 1.6 there):

\begin{theorem}\label{th:dj}
Let $K$ be a convex origin-symmetric body in $\RR^d$. Then there exists a  symmetric, proper,  generalized arithmetic progression ${\bf P} \subset \ZZ^d$, such that $\mbox{rank}({\bf P}) \le d$
and
\begin{equation}\label{eq:incl}
(O(d)^{-3d/2} K) \cap \ZZ^d \subset \text{\rm Image}({\bf P}) \subset K \cap \ZZ^d,
\end{equation}
in addition
\begin{equation}\label{eq:size}
O(d)^{-7d/2} \#K \le \#{\bf P}.
\end{equation}
\end{theorem}

Now we are ready to state and prove our first estimate  in Koldobsky's question and prove that
for any origin-symmetric convex body $K \subset \RR^d$,
$\dim(K\cap\ZZ^d)=d$, 
\begin{equation}\label{eq:viajohn}
\#K\leq O(d)^{7d/2}  \max_{\vxi\in \SP^{d-1}} \left( \#(K\cap \vxi^\perp)\right) \vol_d(K)^{\frac{1}{d}}.
\end{equation}
 To prove  (\ref{eq:viajohn}) we apply the discrete  John's theorem  to get a symmetric, proper,  generalized arithmetic progression ${\bf P}=(N, \vv, d)$ as in Definition \ref{def:progression}. We note that if $\mbox{rank}({\bf P}) < d$, then
there exists a hyperplane $\vxi^\perp$ such that ${\bf P} \subset \vxi^\perp$ and using (\ref{eq:size}) we get
$$
O(d)^{-7d/2} \#K \le \#({\bf P}) \leq \#(K\cap \vxi^\perp).
$$
By our assumption $\dim (K\cap\ZZ^d)=d$ we have $\vol_d(K)\geq 2^d/d!$
and so  $\vol_d(K)^{\frac{1}{d}}> 2/d$. Thus 
$$
 \#K \leq O(d)^{7d/2}\,\#(K\cap \vxi^\perp)\,\vol_d(K)^{\frac{1}{d}}.
$$
Next we  consider the case $\mbox{rank}({\bf P})=d$.
Without loss of generality, take $N_1\geq N_2\geq \ldots \geq N_d \ge 1$, then define $\vxi^\perp=\spn\{\vv_1, \dots, \vv_{d-1}\}$.
Application of  (\ref{eq:size}) gives  
\begin{align*}
\#K  \leq & O(d)^{7d/2}  \#({\bf P}) \\
\leq &O(d)^{7d/2} \prod_{i=1}^{d} (2N_i+1)\\
= &O(d)^{7d/2} (2N_d+1) \prod_{i=1}^{d-1} (2N_i+1)\\
\leq  &O(d)^{7d/2} \left(\prod_{i=1}^{d} (2N_i+1)\right)^{\frac{1}{d}} \#(K\cap \vxi^\perp).
\end{align*}
Where the last inequality follows from the minimality of $N_d$ and we use  (\ref{eq:incl}) to claim that
$$\#(K\cap \vxi^\perp) \geq \prod_{i=1}^{d-1} (2N_i+1).$$
Now we consider the volume covered by our progression.  Take a fundamental parallelepiped
$$\Pi= \left\{a_1 \vv_1 + \ldots +a_d \vv_d, \mbox{ where } a_i \in [0,1), \mbox{ for all }  i=1, \dots, n \right\}.$$
Let $X=[-N,N-1]\cdot \vv$, we notice that 
$$
K\supset \bigcup\limits_{\vx\in X} (\vx+ \Pi),
$$ 
indeed from $\text{\rm Image}({\bf P})\subset K\cap \ZZ^d$ we get that the vertices of $\vx+\Pi$ belong to  $K \cap \ZZ^d$ for all $\vx\in X$  and thus, by convexity,  $\vx+\Pi\subset K$ for all $\vx\in X$. Next
 $$
 \vol_d(K) \geq \left(\prod_{k=1}^d 2 N_k \right)  \det(\vv_1,\ldots, \vv_d) \geq \prod_{k=1}^d 2 N_k,
 $$
 where the last inequality follows from the fact that $\vv_1,\ldots, \vv_d$ are independent vectors in $\ZZ^d$ and thus
  $\det(\vv_1,\ldots, \vv_d) \geq \det(\ZZ^d)=1$.

Finally,
\begin{align*}\#K \leq   & O(d)^{7d/2} \left(\prod_{i=1}^{d} (2N_i+1)\right)^{\frac{1}{d}} \#(K\cap \vxi^\perp) \\ \leq  & O(d)^{7d/2} \left(\prod_{i=1}^{d} (2N_i)\right)^{\frac{1}{d}} \#(K\cap \vxi^\perp)
\\\leq & O(d)^{7d/2}   \#(K\cap \vxi^\perp) \vol_d(K)^{\frac{1}{d}}.
\end{align*}


%
\section{The case of co-dimensional slices and improved bound on $C(d)$}
The  goal of this section is to improve the  estimate provided in Section 3.  We will need to consider counting points intersecting a body with a  general lattice, and so we will adapt our notation slightly. We refer to \cite{TV2}, \cite{BR} and \cite{HW} for the general facts and introduction on the properties of the cardinality of intersections of convex bodies and a lattice.  Given a lattice $\Gamma$ we will take $\#(K\cap \Gamma)= \text{card}(K\cap \Gamma)$ and, as before, if the lattice is omitted we will take the lattice to be the standard integer lattice of appropriate dimension. We begin with the statement of Minkowski's Second Theorem which is an extension of  Minkowski's First Theorem (Theorem \ref{th:M1} above) and can be found, for example, in \cite{HW} (Theorem 1.2) or \cite{TV2} (Theorem 3.30 pg 135). First we recall the definition of Successive Minima.
\begin{definition}
Let $\Gamma$ be a lattice in $\RR^d$ of rank $k$, and let $K$ be an 
origin-symmetric 
convex body in $\RR^d$.  For $1\leq j\leq k$ define the successive minima to be
$$\lambda_j=\lambda_j(K, \Gamma)= \min\left\{\lambda>0 : \lambda \cdot
  K \text{ contains } j \text{ linearly independent elements of }
  \Gamma \right\}.$$
\end{definition}
Note, that it follows directly from the definition that $\lambda_k \ge
\lambda_{k-1}\ge...\ge\lambda_1$. In addition, the  assumption that $K$ contains $d$ linearly independent lattice points of $\Gamma$ implies that  $\Gamma$ has rank $d$ and that $\lambda_d\leq 1$.  Moreover, according to the definition of the
successive minima there exists a set of linearly independent vectors
from $\Gamma$, $\vv_1,\ldots,\vv_k$, such that $\vv_i$ lies on the
boundary of $\lambda_i \cdot K$ but the interior of $\lambda_i \cdot
K$ does not contain any lattice vectors outside the span of
$\vv_1,\ldots, \vv_{i-1}$.  The  vectors $\vv_1,\dots,\vv_k$ are called a directional
basis, and we note that they may not necessarily form  a basis of $\Gamma$.

\begin{theorem}\label{th:M2} (Minkowski's Second Theorem)
Let $\Gamma$ be a lattice in $\RR^d$ of rank $d$, $K$ be an
origin-symmetric convex body with successive minima $\lambda_i$. Then,
$$
\frac{1}{d!} \prod_{i=1}^{d} \frac{2}{\lambda_i} \leq \frac{\vol(K)}{\det(\Gamma)} \leq \prod_{i=1}^{d} \frac{2}{\lambda_i}.
$$
\end{theorem}

Next we will study the behavior of constant $\LLL_3$ in the case of unconditional convex bodies. A set $K\subset \RR^d$ is said to be unconditional if it is symmetric with respect to any coordinate hyperplane, i.e., $(\pm x_1, \pm x_2, \dots, \pm x_d) \in K$, for any $\vx \in K$ and any choice of $\pm$ signs.
\begin{theorem}\label{th:unconditional}
Let $K\subset\RR^d$ be an unconditional convex body with
$\dim(K\cap\ZZ^d)=d$.  Then
$$
\#K\leq O(d) \max_{i=1,\dots,d} \left( \#(K\cap \ve_i^\perp)\right)\,\vol_d(K)^{\frac{1}{d}},
$$
 where $\ve_1, \dots \ve_d$ are the standard basis vectors in
 $\RR^d$. Moreover, this bound is the best possible. 
\end{theorem}
\beginproof
This result follows from the simple observation that the section of $K$ by a coordinate hyperplane
${\ve_i}^\perp$ is maximal in cardinality among all parallel sections of $K$, i.e.
\begin{equation}\label{eq:uncon}
\#(K \cap (\ve_i^\perp + t \ve_i))\le \#(K \cap \ve_i^\perp), \mbox{ for all }  t\in \RR,  \mbox{ and }  i=1,\dots,d.
\end{equation}
We can see this by considering a point $\vx \in K \cap (\ve_i^\perp + t \ve_i)\cap \ZZ^d$.
Let $\bar{\vx}$ be the reflection of $\vx$ over ${\ve_i}^\perp$, i.e., $\bar{\vx}=(x_1, \dots, -x_i, \dots, x_d)$. Using  unconditionality of $K$, we get that $\bar{\vx}\in K$
and convexity gives us $(\vx+\bar{\vx})/2 \in K \cap \ve_i^\perp$.
Hence, the projection of a point in $K \cap (\ve_i^\perp + t \ve_i)$ is associated to a point in $K \cap \ve_i^\perp$, which explains (\ref{eq:uncon}).

Let  $\{\lambda_i\}_{i=1}^d$ be the successive minima of $K$ with
respect to $\ZZ^d$.  
Using an argument similar to the one above one can show that  that the  vectors $\vv_1, \dots, \vv_d \in \ZZ^d$  associated with $\{\lambda_i\}_{i=1}^d$  may be taken as a rearrangement of $\ve_1,\dots, \ve_d$.
We may assume without loss of generality that $\lambda_d$ corresponds to $\ve_d$. So $\ve_d \in \lambda_d K$ and
$\frac{1}{\lambda_d} \ve_d \in K$. Thus $\#(K\cap L_{\ve_d}) \le 2 \lfloor \frac{1}{\lambda_d}\rfloor +1$, where, as before, $L_{\ve_d}$ is a line containing $\ve_d$ and the origin.  Using (\ref{eq:uncon}), we get
$$
\#K\leq  \left( 2 \left\lfloor \frac{1}{\lambda_d}\right\rfloor +1  \right) \#(K\cap {\ve_d}^\perp).
$$
By assumption we have  $\lambda_d  \le 1$ and, using $\lambda_d \ge \lambda_i$, for all $i=1, \dots, d$, we get
$$
 2 \left\lfloor \frac{1}{\lambda_d}\right\rfloor +1  \le \frac{3}{\lambda_d} \le O(d) \left(\frac{1}{d!} \prod_{i=1}^{d} \frac{2}{\lambda_i}  \right)^{1/d}.
$$
Finally we use Theorem \ref{th:M2} to finish the proof:
$$
\#K\leq O(d)    \#(K\cap {\ve_d}^\perp) \,\vol_d(K)^{\frac{1}{d}}.
$$
The cross-polytope $B_1^d =\conv\{\pm\ve_1,\dots,\pm\ve_d\}$ of
$\vol(B_1^d)=2^d/d!$ shows that the bound is optimal up to
multiplication with  constants. 
\endproof

The idea of the proof of the above theorem follows from the classical Brunn's theorem: the central hyperplane section of a convex origin-symmetric body is maximal in volume among all parallel sections (see \cite{Ga}, \cite{K3}, \cite{RZ}).  One may notice that, in general, it may not be the case that the maximal hyperplane in cardinality for an origin-symmetric convex body passes through the origin.  Indeed, see Figure 1 below, or consider an example of a cross-polytope $B_1^d=\{\vx\in \RR^d: \sum|x_i|\le 1\}$, then $\#(B_1^{d} \cap (1/\sqrt{d},\dots, 1/\sqrt{d})^\perp)=1$  but a face of $B_1^d$ contains $d$ integer points.
\begin{figure}[h!]\label{fig:BrunnCounterexample}
\includegraphics[scale=0.4]{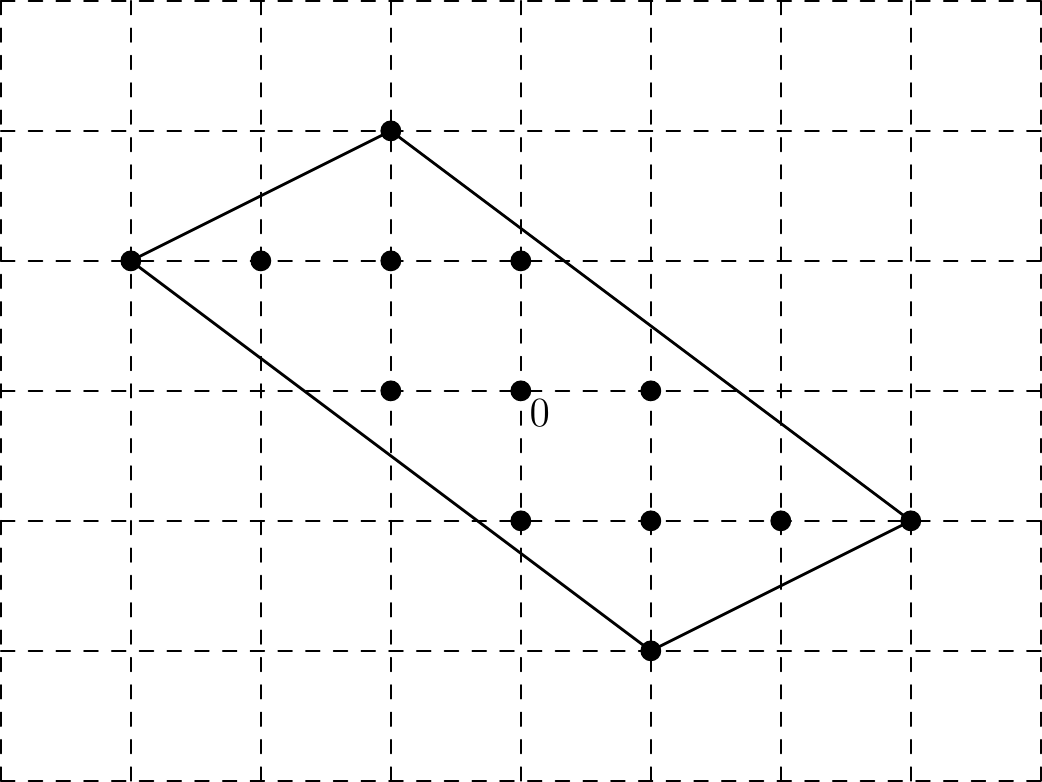}
\caption{Central Section may have less integer points.}
\end{figure}
So we see that there is no equivalent of Brunn's theorem in this setting.

Still, we propose the following analog of  Brunn's theorem in the discrete setting:

\begin{theorem}\label{th:brunn} Consider a convex, origin-symmetric
  body $K\subset \RR^d$ and a lattice $\Gamma \subset \RR^d$ of rank $d$, then
$$
\#(K \cap \vxi^\perp \cap \Gamma) \ge 9^{-(d-1)} \#(K \cap (\vxi^\perp+t \vxi)\cap \Gamma), \mbox{ for all } t \in \RR.  $$
\end{theorem}
\noindent Before proving Theorem \ref{th:brunn} we need to recall two nice packing estimates (see Lemma 3.21, \cite{TV2}):
\begin{lemma}\label{lm:pack}
Let $\Lambda$ be  a lattice in $\RR^d$. If $A \subset \RR^d$ is an arbitrary bounded set and
$P \subset \RR^d$ is a finite non-empty set, then
\begin{equation}\label{eq:3_9}
\#\left(A\cap(\Lambda+P)\right) \le \# \left((A-A)\cap(\Lambda+P-P)\right).
\end{equation}
If $B\subset \RR^d$ is a origin-symmetric convex body, then
\begin{equation}\label{eq:3_10}
(kB) \cap \Lambda\mbox{ can be covered by } (4k+1)^d \mbox{ translates of }B\cap\Lambda.
\end{equation}
\end{lemma}

\noindent{\bf Proof of Theorem \ref{th:brunn}:}
We first recall a standard observation, that the convexity of $K$ gives us
$$
K \cap \vxi^\perp \supset \frac{1}{2}(K \cap (\vxi^\perp+t \vxi))+ \frac{1}{2}(K \cap (\vxi^\perp-t \vxi)).
$$

Let $\Gamma'=\Gamma \cap \vxi^\perp$ and assume that $\Gamma \cap (\vxi^\perp+t \vxi) \not =\emptyset$ (the statement of the theorem is trivial in the other case). Consider  a point  $\vgamma\in \Gamma \cap (\vxi^\perp+t \vxi)$ then
 $$
 \Gamma \cap (\vxi^\perp+t \vxi) = \vgamma+\Gamma'   \mbox{ and }  \Gamma \cap (\vxi^\perp-t \vxi) = -\vgamma+\Gamma'.
 $$
Moreover,
$$
K \cap \vxi^\perp \supset \frac{1}{2}\left(\left[K \cap (\vxi^\perp+t \vxi)\right]-\vgamma\right)+ \frac{1}{2}\left(\left[K \cap (\vxi^\perp-t \vxi)\right]+\vgamma\right).
$$
Thus
 $$
\left(K \cap \vxi^\perp \right) \cap \Gamma' \supset \left[\frac{1}{2}\left(\left[K \cap (\vxi^\perp+t \vxi)\right]-\vgamma\right)+ \frac{1}{2}\left(\left[K \cap (\vxi^\perp-t \vxi)\right]+\vgamma\right)  \right] \cap \Gamma'.
$$
Our goal is to estimate the number of lattice points on the right hand side of the above inclusion. Let
$$
B=\frac{1}{2}\left(\left[K \cap (\vxi^\perp+t \vxi)\right]-\vgamma\right)
$$
then, using the symmetry of $K$, we get
$$
-B=\frac{1}{2}\left(\left[K \cap (\vxi^\perp-t \vxi)\right]+\vgamma\right).
$$
Thus $B-B$ is an origin-symmetric convex body in $\vxi^\perp$. Next we use (\ref{eq:3_10}) from Lemma \ref{lm:pack}  to claim that
$$
\#(2(B-B) \cap \Gamma') \le 9^{d-1} \#((B-B)\cap \Gamma').
$$
Notice that $2(B-B)=2B-2B$ thus we may use (\ref{eq:3_9}) from Lemma \ref{lm:pack} with $P=\{{\boldsymbol 0}\}$, $\vxi^\perp$ associated with $\RR^{d-1}$, and $\Lambda=\Gamma^\prime$ to claim that
$$
\#(2(B-B)\cap \Gamma') =\#((2B-2B)\cap \Gamma' )\ge \# (2B\cap \Gamma')=\#(2B\cap \Gamma).
$$
Thus we proved that
\begin{align*}
\#\bigg(\bigg[\frac{1}{2}(K \cap & (\vxi^\perp+t \vxi) -\vgamma)+ \frac{1}{2}(K \cap (\vxi^\perp-t \vxi)+\vgamma)  \bigg] \cap  \Gamma \bigg)  \\ \ge & 9^{-(d-1)} \# \left( \left[K \cap (\vxi^\perp+t \vxi)-\vgamma  \right]\cap \Gamma \right)
\end{align*}
but 
$$
\#\left(\left[K \cap (\vxi^\perp+t \vxi)-\vgamma  \right] \cap \Gamma \right)=\#\left(\left[K \cap (\vxi^\perp+t \vxi) \right]\cap\Gamma \right).
$$
\endproof

\begin{corollary}\label{cor:Brunn}
Consider a convex, origin-symmetric body $M\subset \RR^n$,  lattice
$\Lambda \subset \RR^n$ and $m$-dimensional lattice subspace $H$,
i.e., it contains $m$ linearly independent  points of $\Lambda$,   then
$$
\#(M \cap H \cap \Lambda) \ge 9^{-m} \#(M \cap (H+\vz)\cap \Lambda), \mbox{ for all } \vz \in \RR^n.  $$
\end{corollary}
\beginproof  Let $\vz\in\RR^n$. Then we may assume $\vz \in \Gamma\setminus\{H\cap\Gamma\}$   and  let  $U$ be the linear space spanned by $H$ and $\vz$.  Then $\dim(U)=m+1$ and the corollary follows  from Theorem \ref{th:brunn} with $U$ associated with $\RR^{m+1}$, $K= M \cap U$, and $\Gamma =\Lambda \cap U$.
\endproof

Let $\G_{\ZZ}(i,d)$ be the set of all $i$-dimensional linear
subspaces containing $i$-linearly independent lattice vectors of $\ZZ^d$, i.e., the
set of all $i$-dimensional lattice hyperplanes.   The next theorem gives a general bound on the number of integer points in co-dimensional slices.

\begin{theorem} 
 Let $K \subset \RR^d$ be an origin-symmetric convex body with $\dim
 (K\cap\ZZ^d)=d$. Then 
\begin{equation} 
 \# K \leq O(1)^d\, d^{d-m} \max\{\# (K\cap H): H\in
 \G_{\ZZ}(m,d)\}\,\, \vol_d(K)^{\frac{d-m}{d}}. 
\end{equation} 
\label{tm:main}
\end{theorem}
Obviously, for $m=d-1$ we obtain the estimate for hyperplane slices 
\begin{equation}
\#K\leq O(1)^d  \max_{\vxi\in \SP^{d-1}} \left( \#(K\cap
  \vxi^\perp)\right)\,\vol_d(K)^{\frac{1}{d}}.
\label{eq:slice_bound} 
\end{equation} 
\medskip

\beginproof
Let  $\{\lambda_i^*\}_{i=1}^d$ be the successive minima of the polar
body $$K^*=\{\vy\in\RR^d: \vy\cdot,\vx\leq 1, \text{ for all }\vx\in
K\}$$ with respect to $\ZZ^d$ and let $\vv_1, \dots \vv_d \in \ZZ^d$ be
the  associated directional basis. These vectors are linearly independent and
$\vv_i\in\lambda_i^*\,K^*$ for all $i$. Thus we have 
\begin{equation}\label{eq:polar}
K\subseteq \{\vx\in\RR^d : |\vv_i \cdot \vx|\leq \lambda_i^*, \,
1\leq i\leq d\}. 
\end{equation}

 Let $U=\spn\{\vv_1,\dots,\vv_{d-m}\}$ and let $\overline H=U^\perp$
 be the orthogonal complement of $U$. Observe that $\overline H\in
 \G_{\ZZ}(m,d)$. Since for $\vz\in\ZZ^d$ we have $\vv_i \cdot
 \vz\in\ZZ$, $1\leq i\leq d$, we also have $\vv_i \cdot (\vz\big
 |U)\in\ZZ$, $1\leq i\leq d-m$, where $\vz\big |U$ is the orthogonal
 projection onto $U$.    In view of \eqref{eq:polar} we obtain  
\begin{equation}
\label{eq:bound_pro} 
   (K\cap\ZZ^d)\big|U \subset \{\vy\in U : \vv_i \cdot \vy\in\ZZ\text{ and }| \vv_i\cdot \vy|\leq \lambda_i^*, \,
1\leq i\leq d-m\}, 
\end{equation}
and thus 
\begin{equation}
  \#((K\cap\ZZ^d)\big|U)\leq \prod_{i=1}^{d-m}\left(2\,\lfloor\lambda_i^*\rfloor +1\right).  
\end{equation}
Due to our assumption that $K$ contains $d$-linearly independent lattice points we have that 
$\lambda_1^*\geq 1$; otherwise  \eqref{eq:polar} implies $\vv_1 \cdot
\vz=0$ for all $\vz\in K\cap\ZZ^d$.  So we conclude by \eqref{eq:bound_pro} 
\begin{equation}
\begin{split}
\# K & \leq \#((K\cap\ZZ^d)|U)\,\max\{\#(K\cap (\vz+\overline H)): \vz\in\ZZ^d \} \\
     & \leq \max\{\#(K\cap (\vz+\overline H)): \vz\in\ZZ^d \}\,3^{d-m}\prod_{i=1}^{d-m}\lambda_i^*  \\
     &\leq 3^{d-m}\,O(1)^m \#(K\cap \overline H) \,\prod_{i=1}^{d-m}\lambda_i^* 
      \leq O(1)^d\,\#(K\cap \overline H) \,\prod_{i=1}^{d-m}\lambda_i^* . 
\end{split}
\label{eq:bound1}
\end{equation} 
Here the last step follows from Corollary \ref{cor:Brunn}, the  co-dimensional version of the discrete Brunn's theorem.

Next Minkowski's Second theorem (Theorem \ref{th:M2})  gives the upper bound 
\begin{equation}
              \lambda_1^*\cdot\ldots\cdot\lambda_d^* \vol_d(K^*)\leq 2^d
\end{equation}
and so we find 
\begin{equation}
\left(\prod_{i=1}^{d-m}\lambda_i^*\right)^d\vol_d(K^*)^{d-m}\leq \left(\prod_{i=1}^{d}\lambda_i^*\right)^{d-m}\vol_d(K^*)^{d-m}\leq 2^{d(d-m)}.
\end{equation}
Hence 
\begin{equation}
\prod_{i=1}^{d-m}\lambda_i^* \leq 2^{d-m}\, \vol_d(K^*)^{\frac{m-d}{d}}. 
\label{eq:mink} 
\end{equation}
By the Bourgain-Milman inequality (isomorphic version of reverse
Santal\'{o} inequality see \cite{BM, GPV, Na, Ku} or \cite{RZ}),  there
exists  an absolute constant $c>0$ with 
 $$
 c^d \frac{4^d}{d!} \leq \vol_d(K)\vol_d(K^*)
 $$
and so we get 
\begin{equation}
        \vol_d(K^*)^{\frac{m-d}{d}} \leq O(d)^{d-m}\vol_d(K)^{\frac{d-m}{d}}.
\end{equation} 
Thus together with \eqref{eq:mink} and \eqref{eq:bound1} we obtain 
\begin{equation}
\# K \leq O(1)^d\,d^{d-m}\, \max\{\# (K\cap H): H\in
 \G_{\ZZ}(m,d)\}\, \vol_d(K)^{\frac{d-m}{d}}.
\end{equation}
\endproof

\begin{remark} We notice that the methods used in Section 3,
  i.e. computation via discrete version of the John theorem (Theorem
  \ref{th:dj} from above), can also be used to provide a bound for
  general co-dimensional  sections. But such computation gives the estimate of order $O(d)^{7d/2}$ which is worse than the one in the above theorem.
\end{remark}

\begin{remark} Observe that Theorem \ref{tm:main} can be restated for  an arbitrary  $d$-dimensional lattice $\Lambda$: Let $\Lambda$ be a lattice in $\RR^d$ and $K \subset \RR^d$ be an origin-symmetric convex body with $\dim
 (K\cap \Lambda)=d$. Then 
\begin{equation} 
 \# K \leq O(1)^d\, d^{d-m} \max\{\# (K\cap H): H\in
 \G_{\Lambda}(m,d)\}\,\, \left(\frac{\vol_d(K)}{\det(\Lambda)}\right)^{\frac{d-m}{d}}. 
\end{equation} 
\end{remark}

We also notice that the methods used in proofs of the Theorem \ref{th:unconditional}  and Theorem \ref{tm:main} can be used to provide an estimate for the co-dimensional slices of an unconditional convex body:

\begin{theorem} 
 Let $K \subset \RR^d$ be an unconditional convex body with $\dim
 (K\cap\ZZ^d)=d$. Then 
\begin{equation} 
 \# K \leq  O(d)^{d-m} \max\{\# (K\cap H): H\in
 \G_{\ZZ}(m,d)\}\,\, \vol_d(K)^{\frac{d-m}{d}}. 
\end{equation} 
\label{tm:unconmain}
\end{theorem} 
\beginproof  First we notice that if $K$ is an unconditional body and $H$ is a coordinate subspace of dimension $m$ (i.e. it is spanned by $m$ coordinate vectors) with $
K \cap (H + \vz) \not = \emptyset,
$
 then $
K \cap (H + \vz) 
$ must be an unconditional convex body in  $(H + \vz)$. Thus, using this property together  with  the proof of Theorem \ref{th:unconditional} we get  that  for any unconditional body $K$ and for any coordinate subspace $H$
$$
\#(K \cap H \cap \ZZ^d) \ge  \#(K \cap (H+\vz)\cap \ZZ^d), \mbox{ for all } \vz \in \RR^d.
$$
Next we follow the steps of the proof of Theorem \ref{tm:main} and similarly to   (\ref{eq:bound1})
 get 
$$
\# K \leq 3^{d-m} \#(K\cap \overline H) \,\prod_{i=1}^{d-m}\lambda_i^*. 
$$
Finally, we finish the proof using Minkoswki's Second theorem and the Bourgain-Milman inequality.

\endproof

\begin{remark} We also would like to test our estimates against two classical examples 
\begin{enumerate}
\item[{\bf(A)}] For the cube $B_\infty^d = \{\vx\in\RR^d : |\vx|_\infty\leq 1\}$ we have 
$\# B_\infty^d =3^d$,\\
 $\max\{\# (B_\infty^d\cap H): H\in \G_{\ZZ}(m,d)\}=3^m $, and $\vol_d(B_\infty^d)^{\frac{d-m}{d}} = 2^{d-m}$.

\item[{\bf (B)}] For the cross polytope $B_1^d = \{\vx\in\RR^d : |\vx|_1\leq 1\}$ we have 
$\# B_1^d =2\,d+1$,\\
 $\max\{\# (B_1^d\cap H): H\in \G_{\ZZ}(m,d)\}=2\,m+1$, and $\vol_d(B_1^d)^{\frac{d-m}{d}} \sim \frac{c^{d-m}}{d^{d-m}}$.
\end{enumerate} 
These examples show that we expect our constant to grow exponentially in the case of higher co-dimensional slices, though we do not expect our current estimates to be sharp.
\end{remark}

We finish this section with a remark  about the relationship between the constant in the original slicing inequality, $\LLL_1$, and the constant in the discrete version, $\LLL_3$.  Using the general idea from \cite{GGroZ} and Gauss's Lemma on the intersection of a large convex body with a lattice we will show that $\LLL_1 \le \LLL_3$.

Consider a convex symmetric body $K$ and let $\LLL_1(K)>0$ be such that
$$\vol_d(K)^{\frac {d-1}{d}} = \LLL_1(K) \max_{\vxi \in \SP^{d-1}} \vol_{d-1}(K\cap \vxi^\bot).$$
Thus $$\LLL_1=\max \{\LLL_1(K): K \subset \RR^d,  K \mbox{ is convex, origin-symmetric body, } d \ge 1\}.$$ Then  $$\vol_d(K)^{\frac {d-1}{d}} \geq \LLL_1(K) \vol_{d-1}(K\cap \vxi^\bot), \,\, \forall \vxi\in \SP^{d-1}.$$
Our goal is to study a section of $K$ with a maximal number of points from $\ZZ^d$, if $K\cap \vxi^\perp$ is such a section, then, without loss of generality,  we may assume that $\ZZ^d \cap \vxi^\perp$ is a lattice of a full rank $d-1$. Indeed, if $\ZZ^d \cap \vxi^\perp$ has a rank less then $d-1$ we may rotate $\vxi$ to catch $d-1$ linearly independent vectors in $\vxi^\perp$, without decreasing the number of integer points in $K\cap \vxi^\perp$. Now, we may use Gauss's Lemma (see for example Lemma 3.22 in \cite{TV2}) to claim that for $r$ large enough we have
$$\#(rK) = r^d \vol_d(K) + O\left(r^{d-1}\right) \text{ and}$$
$$\#(rK \cap \vxi^\perp)=\frac{r^{d-1}\vol_{d-1}(K)}{\det\left(\ZZ^d \cap \vxi^\perp\right)} + O\left(r^{d-2}\right).$$
Which we can rearrange to get the following two equations
$$\vol_{d}\left(K\right) = \frac{1}{r^{d}} \#\left(rK \right) + O\left(\frac{1}{r}\right) \text{ and}$$
$$\vol_{d-1}\left(K\cap \vxi^\perp\right) = \frac{1}{r^{d-1}} \#\left(rK \cap \vxi^\perp\right) \det\left(\ZZ^d\cap \vxi^\perp\right) + O\left(\frac{1}{r}\right).$$
Next, using that  $\det\left(\ZZ^d \cap \vxi^\perp \right)\geq 1$ and  $\vol_d(K) \geq \LLL_1(K) \vol\left(K \cap \vxi^\perp\right) \vol_d^{\frac{1}{d}}(K)$ we get
\begin{align*}
\frac{1}{r^{d}} \#\left(rK \right)  \geq& \LLL_1(K) \left(\frac{1}{r^{d-1}}\right) \#\left(rK \cap \vxi^\perp\right) \det\left(\ZZ^d\cap \vxi^\perp\right)\vol_d^{\frac{1}{d}}(K) + O\left(\frac{1}{r}\right)\\
\geq& \LLL_1(K) \left(\frac{1}{r^{d}}\right)\#\left(rK \cap \vxi^\perp\right)\vol_d^{\frac{1}{d}}(rK)+ O\left(\frac{1}{r}\right).
\end{align*}
Then for $\epsilon>0$ there is a sufficiently large $r_0$ such that for all $r>r_0$
$$\#(rK)\geq \left(\LLL_1(K)-\epsilon\right)\#\left(rK \cap \vxi^\perp\right)\vol_d^{\frac{1}{d}}(rK).$$
So then if $\#(rK) \leq \LLL_3 \max_{\vxi\in \SP^{d-1}} \# \left(rK \cap \vxi^\perp\right) \vol_d^{\frac{1}{d}}(rK)$ we have that $\LLL_1(K)-\epsilon \leq \LLL_3$ for all $d$, $\epsilon$, and bodies $K$.
Which leads us to conclude that $\LLL_1\leq \LLL_3$.



\end{document}